\documentclass[11pt]{amsart}
\usepackage[utf8]{inputenc}
\usepackage{amsfonts}
\usepackage{amssymb}
\usepackage{amsmath}
\usepackage{mathtools}
\usepackage{xcolor}
\usepackage{amsthm}
\usepackage{pdflscape}
\usepackage{pgfplots}
\usepackage{mathrsfs}
\usepackage{hyperref}
\usepackage{comment} 
\usepackage[pagewise]{lineno}

\newtheorem{thm}{Theorem}[section]
\newtheorem{cor}[thm]{Corollary}
\newtheorem{prop}[thm]{Proposition}
\newtheorem*{prop1}{Proposition}
\newtheorem*{thm1}{Theorem}
\newtheorem{lemma}[thm]{Lemma}

\newtheorem{fact}[thm]{Fact}

\theoremstyle{definition}
\newtheorem{defn}[thm]{Definition}

\newtheorem{notn}[thm]{Notation}
\newtheorem{remark}[thm]{Remark}

\newcommand{\gen}[1]{\left\langle#1\right\rangle}
\theoremstyle{remark}
\newtheorem{rem}[thm]{Remark}

\theoremstyle{plain}

\newcommand{\Mod}[1]{\;(\text{mod}\;#1)}

\begin{document}

\title{On transformations in the Second Painlev\'e family}
\author{Christine Eagles}
\address{Christine Eagles\\
University of Illinois at Chicago\\
Department of Mathematics, Statistics and Computer Science\\
Science and Engineering Offices\\
Chicago, IL, 60607\\
United States}
\email{ceagles@uic.edu}

\author{Joel Nagloo}
\address{Joel Nagloo\\
University of Illinois at Chicago\\
Department of Mathematics, Statistics and Computer Science\\
Science and Engineering Offices\\
Chicago, IL, 60607\\
United States}
\email{jnagloo@uic.edu}

\keywords{}
\subjclass[2020]{}

\thanks{J. Nagloo is partially supported by NSF grant DMS-2348885. C. Eagles is partially supported by the Natural Sciences and Engineering Research Council of Canada (NSERC), [CPRA - 609532 - 2026].}

\begin{abstract} 
\sloppy We complete the classification of algebraic relations among solutions of equations in the second Painlev\'e family $P_{II}(\alpha): y''=2y^3+ty+\alpha$, where $\alpha\in\mathbb{C}$. More precisely, we prove that if $y_1,\ldots,y_n$ are generic solutions of $P_{II}(\alpha_1),\ldots,P_{II}(\alpha_n)$, respectively, and $tr.deg_{\mathbb{C}(t)}\mathbb{C}(t)(y_1,y_1',\ldots,y_n,y_n')<2n$, then for some $i<j$, we have that $\alpha_i\pm\alpha_j\in\mathbb{Z}$, and $y_i=T(y_j)$ for a classical B\"acklund transformation $T$. In particular, the classical B\"acklund transformations are the only B\"acklund transformations in the second Painlev\'e family.
\end{abstract}

\maketitle

\tableofcontents

\section{Introduction}

Ever since the seminal work of Painlev\'e in the early 1900s on the existence of the six equations that now bear his name, mathematicians have been fascinated and curious about the transcendental nature of the solutions to these equations, namely, the Painlev\'e transcendents. Indeed, one of the earliest claims of Painlev\'e about these equations (cf. \cite{Pain}) was that generically, their solutions are truly new transcendental functions: they cannot be expressed in terms of elementary or {\em classical functions}.

In the early 1980s, this claim of Painlev\'e was taken up by the Japanese school of differential algebra and by the end of the 1990s, the work of Nishioka \cite{Nishioka1}, Okamoto-Noumi \cite{NoumiOkamoto}, Umemura-Watanabe \cite{Umemura2}, \cite{Umemura3}, and Watanabe \cite{Watanabe}, \cite{Watanabe1},  had given a complete classification of the classical solutions to each of the Painlev\'e equations. In so doing, they established that indeed, except for special values of the parameters, the equations are irreducible in the sense of not having classical solutions. Later, Nagloo-Pillay \cite{NagPil} discovered that what the Japanese school had actually established is that generically, the sets defined by the Painlev\'e equations are strongly minimal; a fundamental property in model theory, and central to the modern interaction between model theory and differential algebra.

\sloppy One of the remarkable consequences of the trichotomy theorem of Hrushovski-Sokolovic \cite{Hrushovski-Sokolovic} is that the possible algebraic relations among solutions of a strongly minimal equation are subject to very strong structural constraints. For an order two equation, such as a Painlev\'e equation, either they are closely related to the solution set of a very special kind of ODE on an elliptic curve (the Painlev\'e-Picard is the most well-known example) or they are geometrically trivial, namely, all algebraic relations, if any, arise between {\em pairs} of solutions and their derivatives. This naturally raised a new question: having established that the Painlevé transcendents are genuinely new functions, what algebraic relations can exist between them?

Nishioka \cite{Nishioka} completely answered this question for the first Painlev\'e equation. He proved that if $y_1,\ldots,y_n$ are distinct solutions of $P_{I}:= y''=6y^2+t$ then $tr.deg_{\mathbb{C}(t)}\mathbb{C}(t)(y_1,y'_1,\ldots,y_n,y'_n)=2n$, that is, the solutions and their derivatives are algebraically independent over $\mathbb{C}(t)$. For the other Painlev\'e equations $P_{II}-P_{VI}$, the situation is much more delicate. It is well-known that there are the so-called classical B\"acklund transformations that take solutions from one Painlev\'e equation to another within the same family (possibly with the same parameters).
For example, if $z$ is a solution of $P_{II}(\alpha)$, then it is not hard to see that
\begin{eqnarray*}
S(z)&=&-z\\
T_{+}(z)&=&-z-\frac{\alpha+1/2}{z'+z^2+t/2}\\
T_{-}(z)&=&-z+\frac{\alpha-1/2}{z'-z^2-t/2}
\end{eqnarray*}
are solutions of $P_{II}(-\alpha)$, $P_{II}(\alpha+1)$ and $P_{II}(\alpha-1)$, respectively. These transformations generate a group $G$ isomorphic to the affine Weyl group of type $\tilde{A}_1$, and analogous groups of transformations exist for the other Painlev\'e equations $P_{III}-P_{VI}$. Such B\"acklund transformations played a crucial role in the classification of algebraic and classical solutions carried out by the Japanese school mentioned above.

The existence of these classical B\"acklund transformations has a natural explanation in the theory of isomonodromic deformations. Recall that each Painlev\'e equation arises as the compatibility condition of a pair of linear differential systems, known as a Lax pair. It follows that the Schlesinger transformations – special rational gauge transformations of these linear systems – preserve their monodromy data while shifting the parameters, thereby inducing B\"acklund transformations of the corresponding Painlev\'e equations.

So for $P_{II}-P_{VI}$, the natural question is: are these classical transformations responsible for every algebraic relation between Painlev\'e transcendents\footnote{For $P_{VI}$ one also has to account for the Painlev\'e-Picard solutions.}? In this paper, we answer this question in the affirmative for the second Painlev\'e equation $P_{II}(\alpha)$. Let us recall the history of progress made on this problem for $P_{II}$.

Nagloo-Pillay \cite{NagPil} showed that if $\alpha$ is transcendental, then the set defined by $P_{II}(\alpha)$ in a differentially closed field of characteristic zero is geometrically trivial. Namely, they proved that if $y_1,\ldots,y_n$ are distinct solutions of $P_{II}(\alpha)$ and $tr.deg_{\mathbb{C}(t)}\mathbb{C}(t)(y_1,y'_1,\ldots,y_n,y'_n)<2n$, then already for some $i<j$ we have $tr.deg_{\mathbb{C}(t)}\mathbb{C}(t)(y_i,y'_i,y_j,y'_j)=2$. Later, Nagloo \cite{Nagloo} extended this result to all $\alpha\not\in\frac{1}{2}+\mathbb{Z}$, that is, to the case where $P_{II}(\alpha)$ is strongly minimal and, in particular, has no order 1 subvariety over $\mathbb{C}(t)$. Building on this work, Nagloo-Pillay \cite{NagPil} showed that for $\alpha$ transcendental, no such pairwise relations can exist; in other words, $tr.deg_{\mathbb{C}(t)}\mathbb{C}(t)(y_1,y'_1,\ldots,y_n,y'_n)=2n$ for any collection of distinct solutions $y_1,\ldots,y_n$ of $P_{II}(\alpha)$. This completely settled the question of algebraic relations between the solutions of a fixed generic $P_{II}(\alpha)$, but did not address the problem for algebraic $\alpha\in\mathbb{Q}^{alg}$.

The question of algebraic relations between the solutions of distinct equations $P_{II}(\alpha_1),\ldots,P_{II}(\alpha_n)$ was taken up by Nagloo \cite{Nagloo2}, and Freitag-Nagloo \cite{FreitagNagloo}. In the cases where geometric triviality holds, the problem is reduced to classifying the algebraic relations between pairs of equations $P_{II}(\alpha)$ and $P_{II}(\beta)$. When such relations exist, we say that $P_{II}(\alpha)$ and $P_{II}(\beta)$ are nonorthogonal. Freitag and Nagloo proved the deepest result in the case of transcendental parameters\footnote{The paper \cite{FreitagNagloo} claims a stronger result. However, as we will point out in Remark \ref{mistake}, there is a mistake in the proof that $P_{II}(\alpha)$ is geometrically trivial when $\alpha\in\frac{1}{2}+\mathbb{Z}$.}. 

\begin{prop1}[Freitag-Nagloo \cite{FreitagNagloo}]
Suppose that $\alpha, \beta \in \mathbb{C}$ with $\alpha,\beta$ transcendental. Then $P_{II} ( \alpha )$ is nonorthogonal to $P_{II} (\beta)$ if and only if $\beta - \alpha \in \mathbb{Z}$ or $\beta + \alpha \in \mathbb{Z}$. In the case that $\beta - \alpha \in \mathbb{Z}$ or $\beta + \alpha \in \mathbb{Z}$, the only algebraic relations between solutions of $P_{II} ( \alpha)$ and $P_{II} (\beta)$ are induced by the classical B\"acklund transformations. 
\end{prop1}
The techniques used by Freitag and Nagloo rely crucially on the assumption that $\alpha$ is a transcendental constant. In particular, they do not extend to algebraic parameters and it was clear that a completely different approach needed to be developed to tackle those cases. On the other hand, one could also look at the work of Nishioka \cite{Nishioka} on $P_{I}$ and try to adapt some of the calculations for $P_{II}$. The first step in this direction was already taken by Nagloo \cite{Nagloo}, who showed that any solution $y_1$ of $P_{II}(\alpha)$, and $y_2$ of $P_{II}(\beta)$, with $\alpha,\beta\not\in\frac{1}{2}+\mathbb{Z}$, if $tr.deg_{\mathbb{C}(t)}\mathbb{C}(t)(y_1,y'_1,y_2,y'_2)=2$ then $y_2=f(y_1,y_1')$ for some differential rational function over $\mathbb{C}(t)^{alg}$. In this paper we also tackle the second aspect of Nishioka's approach and show that any such differential rational function over $\mathbb{C}(t)^{alg}$ has a Nishioka normal form
\[
f(Y,Y')=\varepsilon Y+b+\frac{q(Y,Y')'}{q(Y,Y')},
\]
where $\varepsilon=\pm 1$, $b\in\mathbb{C}(t)^{alg}$ and $q\in\mathbb{C}(t)^{alg}(Y,Y')$. The fact that we are able to control the shape of the differential rational function will play a crucial role in this paper. Indeed, it does not take much to see that the transformations $T_+$ and $T_-$ can be put in Nishioka normal form. The question then becomes: how can we reveal the constraints on the coefficients of such a differential rational function?

To approach this question, we turn to the formal Puiseux series solutions of $P_{II}(\alpha)$ at $t=\infty$ found by ChatGPT. We show that, for every $\alpha\in\mathbb{C}$, there are precisely three such solutions. Now let $\hat y$ be one of these solutions. Since the field of formal Puiseux series is algebraically closed, and differential rational functions (over $\mathbb{C}(t)^{alg}$) of Puiseux series are again Puiseux series, the solution $f(\hat y)$ is forced to be one of the three Puiseux series solutions of $P_{II}(\beta)$. In the case $\beta=\alpha$, this already gives, together with geometric triviality, that $P_{II}(\alpha)$ is $\omega$-categorical: if $y$ is a generic solution, then $\mathbb{C}(t)(y,y')^{alg}$ only contains finitely many other solutions of $P_{II}(\alpha)$.

Using the Nishioka normal form of $f$, we are able to predict how $f$ acts on the set of formal Puiseux series solutions of $P_{II}(\alpha)$ at $t=\infty$. This reduces the general problem to the case when $\alpha,\beta\not\in\frac{1}{2}\mathbb{Z}$ and a differential rational bijection $f:P_{II}(\alpha)\rightarrow P_{II}(\beta)$ over $\mathbb{C}(t)$ exists. In that case, we use a strong form of the Painlev\'e property, carefully chosen branches $\hat y_{\alpha}$ and $\hat y_{\beta}$ of $P_{II}(\alpha)$ and $P_{II}(\beta)$, to compare the coefficients of the series in the equation $\hat y_{\beta}=\varepsilon\hat y_{\alpha}+b+\frac{q'}{q}$. This allows us to extract that $\beta\pm \alpha \in \mathbb{Z}$ and forces $f$ to be one of the classical B\"acklund transformation. We get

\begin{thm1}
If $y_1,\ldots,y_n$ are generic solutions of $P_{II}(\alpha_1),\ldots,P_{II}(\alpha_n)$, respectively, and $$tr.deg_{\mathbb{C}(t)}\mathbb{C}(t)(y_1,y_1',\ldots,y_n,y_n')<2n,$$ then for some $i<j$, we have that $\alpha_i\pm\alpha_j\in\mathbb{Z}$, and $y_i=T(y_j)$ for a classical B\"acklund transformation $T$. In particular, the classical B\"acklund transformations are the only B\"acklund transformations in the second Painlev\'e family.
\end{thm1}

The paper is organized as follows. In Section 2 we recall the necessary facts from differential algebra and model theory needed in the paper. In Section 3 we use the Gambier construction to establish geometric triviality at the half-integer parameters. In Section 4 we classify the formal Puiseux series solutions at infinity and prove $\omega$-categoricity for all $\alpha\in\mathbb{C}$. Section 5 deals with transformations within a fixed equation, and finally in Section 6 we complete the classification of transformations between different equations.

\subsection*{AI disclosure} ChatGPT (GPT-5.6 and GPT-6, OpenAI) was used during the research leading to this paper. Building on a known algebraic framework for transformations of the second Painlev\'e equation (strong minimality, geometric triviality, Nishioka normal form, among others), the second author directed ChatGPT through a series of mathematical prompts to investigate the local solutions of the second Painlev\'e equation. ChatGPT helped compute the formal Puiseux solutions at infinity and develop their expansions, and pointed to their relevance to the classification problem, including the analysis of the associated Riccati equations. The authors developed and verified the final mathematical arguments, and they take full responsibility for the contents.



\section{Preliminaries}\label{prelim}

When looking at solutions to systems of differential equations, it will often be convenient to work in $(\mathcal{U},D)$, a universal differential field in the sense of Kolchin. Such a universal domain is also a saturated model of the theory $\mathrm{DCF}_0$ in the language $L=\{0,1,+,-,\times,D\}$. We give some necessary details on differential fields here and refer the reader to \cite{Marker} for a more detailed and comprehensive view. 

First, $(\mathcal{U},D)$ is differential field of characteristic zero. So $\mathcal{U}$ is a field of characteristic zero and $D:\mathcal{U}\to\mathcal{U}$ is a derivation meaning that $D$ is an additive homomorphism and satisfies the Leibniz rule $D(xy)=D(x)y+xD(y)$. By a differential subfield $(K,D')$ of $\mathcal{U}$, we mean a subfield $K$ of $\mathcal{U}$, such that $D|_K=D'.$ We will typically denote the derivation $D(a)=a'$. 

For a differential subfield $K$ of $\mathcal{U}$, we define the ring of differential polynomials over $K$ in the variables $Y_1,....,Y_n$ as $K\{Y_1,...,Y_n\}=K[Y_1,...,Y_n,Y_1',...,Y_n',Y_1^{(2)},...,Y_n^{(2)},...].$ We equip this ring with the natural derivation. This derivation also extends naturally to its function field denoted $K\langle Y_1,...,Y_n\rangle$. By $P^D$, we mean the differential rational obtain by applying the derivation on $K$ to the coefficients of $P$. 

We also require that $\mathcal{U}$ is \emph{differentially closed} meaning that if a finite system of differential polynomials in $\mathcal{U}\{x\}$ has a solution in some differential field extension of $\mathcal{U}$, then it has a solution in $\mathcal{U}$.

Finally, we assume that $\mathcal{U}$ is saturated of cardinality of the continuum and so we identify its field of constants $\mathcal{C}=\{a\in\mathcal{U}\mid a'=0\}$ with the field of complex numbers, $\mathbb{C}$. Furthermore, we may also assume that any finitely generated differential field extension of $\mathbb{C}$ embeds into $\mathcal{U}$. We fix some $t\in\mathcal{U}$ such that $t'=1$. Typically, we will work over $\mathbb{C}(t)$ or $\mathbb{C}(t)^{alg}$.

By a \emph{definable set}, as a consequence of quantifier elimination in $\mathrm{DCF}_0$, we mean a finite Boolean combination of solutions to differential polynomial equations (or Kolchin closed sets). We say $X$ is $K$-definable if the differential polynomial equations defining $X$ have coefficients in $K$. By a $K$-definable function $f:X\to Y$, we mean that the graph of $f$ is $K$-definable. 

Let $K$ be a differential field extension of $\mathbb{C}(t)$. By a \emph{generic solution} of $P_{II}(\alpha)$ over $K$, we mean some $y\in\mathcal{U}$ such that $$P_{II}(\alpha):=\; y''=2y^3+ty+\alpha$$ and $\mathrm{tr.deg}_KK(y,y')=2.$ We will often say that $y$ is a generic solution of $P_{II}(\alpha)$ when $K$ is fixed. Note that thanks to the work of Murata \cite{Murata} and others in the Japanese school of differential algebra, we know exactly when non $\mathbb{C}(t)$-generic solutions exist.

\begin{fact} Let $\alpha\in\mathbb{C}$.
\begin{enumerate}
\item $P_{II}(\alpha)$ has a solution in $\mathbb{C}(t)^{alg}$ if and only if $\alpha\in\mathbb{Z}$. Furthermore, the solution when it exists is unique (and denoted by $r_\alpha$).
\item $P_{II}(\alpha)$ has a solution $y$ with $tr.deg_{\mathbb{C}(t)}\mathbb{C}(t)(y,y')=1$ if and only if $\alpha\in\frac{1}{2}+\mathbb{Z}$.
\end{enumerate}
\end{fact}

If $\alpha=0$, then $r_0:=0$ is the unique rational solution of $P_{II}(0)$. Applying the B\"acklund transformations $T^+$ and $T_-$ consecutively, yield the unique rational solutions $r_{\alpha}$. For example, $r_1=-\frac{1}{t}$. If on the other hand $\alpha=\frac{1}{2}$, it is easy to see that every solution of the Riccati equation $$R_{\frac{1}{2}}:=\;y'=y^2+\frac{t}{2}$$ is also a solution of $P_{II}(\frac{1}{2})$. Using the B\"acklund transformations, we obtain  $R_\alpha$, the {\em exceptional subvariety} for every $\alpha\in\frac{1}{2}+\mathbb{Z}$.

Throughout $X_{II}(\alpha)$ will denote the set defined in $\mathcal{U}$ of the second Painlev\'e equation in the following sense: write
\[
P^*_{II}(\alpha):=
\begin{cases}
P_{II}(\alpha), & \alpha\notin \frac12\mathbb Z,\\[4pt]
P_{II}(\alpha)\setminus\{r_\alpha\}, & \alpha\in\mathbb Z,\\[4pt]
P_{II}(\alpha)\setminus R_\alpha, & \alpha\in\frac12+\mathbb Z,
\end{cases}
\]
where $r_\alpha$ is the unique rational solution and $R_\alpha$ is the exceptional subvariety. Then $X_{II}(\alpha):=\{y\in\mathcal{U}:y\text{ is a solution of }P^*_{II}(\alpha)\}$.

We now recall the known results about the model-theoretic approach to studying the Painlev\'e equations. This will allow us to define the various notions from model theory that will {\em again} play a role in the proof of our main result. First, the connection between the work of the Japanese school described in the introduction and model theory came via the notion of strong minimality.
\begin{defn}
 We say the definable set $X\subset\mathcal{U}$ is \emph{strongly minimal} if any definable subset of $X$ is either finite or cofinite.
\end{defn}

As detailed in \cite[Section 2]{NagPil}, for concrete differential equations, strong minimality has other equivalent characterizations. The following is all we will need 
\begin{fact}\label{factSMConditionJ}
Let $X$ be the set of solutions of
\[
y''=f(t,y,y'),
\]
where $f\in\mathbb C(t)[Y,Y']$. Then the following are equivalent:
\begin{enumerate}
\item $X$ is strongly minimal.
\item For any differential field extension $K$ of $\mathbb C(t)$ and any $y\in X$, either $y\in K^{\mathrm{alg}}$ or $y$ is generic over $K$.
\item The equation satisfies Umemura's Condition $(J)$: for any differential field extension $K$ of $\mathbb C(t)$, if a nonzero $p\in K[y,y']$ divides its derivative $p'$, then $p\in K$.
\end{enumerate}
\end{fact}
Notice that the exceptional subvarieties give the failure of condition $(J)$. Indeed for $\alpha=\frac{1}{2}$, let $y\in X_{II}(\frac{1}{2})$ be any solution (even generic over $\mathbb{C}(t)$) and consider the polynomial $p(y,y')=y'-y^2-\frac{t}{2}$. It is not hard to see that 
\[
\begin{aligned}
p'
&=y''-2yy'-\frac{1}{2}\\
&=2y^3+ty+\frac{1}{2}-2yy'-\frac{1}{2}\\
&=-2y(y'-y^2-\frac{t}{2})\\
&=-2yp,
\end{aligned}
\]
So for a definable set $X$ (or the equation) as in Fact \ref{factSMConditionJ}, failure of Condition (J) over a field $K$ is witnessed by an order 1 equation over $K$ that defines a definable subset. The work of Okamoto-Noumi \cite{NoumiOkamoto} and Umemura-Watanabe \cite{Umemura2} completely tackled the problem for $P_{II}(\alpha)$:
\begin{fact}\label{SM-Invariant}Let $\alpha\in\mathbb{C}$.
\begin{enumerate}
\item $P_{II}(\alpha)$ satisfies Umemura's Condition (J) if and only if $\alpha\not\in\frac{1}{2}+\mathbb{Z}$.
\item If $\alpha\in\frac{1}{2}+\mathbb{Z}$, then $R_{\alpha}$ is the unique definable subset defined by an order 1 equation over any differential field extension $K$ of $\mathbb{C}(t)$.
\end{enumerate}
\end{fact}
Consequently – and this was the observations of Nagloo-Pillay \cite{NagPil} and Freitag-Nagloo \cite{FreitagNagloo} – we get that:
\begin{fact}\label{strongminimality}
$X_{II}(\alpha)$ is strongly minimal for all $\alpha\in\mathbb{C}$. 
\end{fact}

{
\begin{remark}
The equivalence between conditions (1) and (2) in Fact \ref{factSMConditionJ} holds for any definable set. However, for general definable sets, strong minimality is not equivalent to Umemera's condition (J). We note that when $\alpha\in\frac{1}{2}\mathbb{Z}$, then $X_{II}(\alpha)$ is no longer the set of solutions to an equation of the form $y''=f(t,y,y')$ but rather its set of generic solutions.
\end{remark}
}

Once strong minimality is established, classifying the fine structure of $X_{II}(\alpha)$ according to the Trichotomy Theorem of Hrushvoski and Sokolovic is the next natural thing to do. We have the following.

\begin{fact}\label{trivialityP2} Let $\alpha\in\mathbb{C}$.
\begin{enumerate}
\item If $\alpha\not\in\frac{1}{2}+\mathbb{Z}$, then $X_{II}(\alpha)$ is  \em{geometrically trivial}: if $K$ is any differential field extending $\mathbb{C}(t)$ and $y_1,...,y_n\in X_{II}(\alpha)\setminus K^{alg}$, then if $\mathrm{trdeg}_KK(y_1,y_1',....,y_n,y'_n)<2n$, there is some $i<j$ such that $\mathrm{trdeg}_KK(y_i,y_i',y_j,y_j')=2.$
\item If $\alpha\not\in\mathbb{Q}^{alg}$, then $X_{II}(\alpha)$ is \em{strictly disintegrated}: if $K$ is any differential field extending $\mathbb{C}(t)$ and $y_1,...,y_n\in X_{II}(\alpha)\setminus K^{alg}$ are distinct, then $\mathrm{trdeg}_KK(y_1,y_1',....,y_n,y'_n)=2n$.
\end{enumerate}
\end{fact}
As we will soon point out, geometric triviality gives that the field $K$, in the conclusion of both (1) and (2), can be replaced by $\mathbb{C}(t)$. As discussed in the introduction, the proof of Fact \ref{trivialityP2} was first obtained in the case of $\alpha\not\in\mathbb{Q}^{alg}$ by Nagloo-Pillay \cite{NagPil1}. Nagloo \cite{Nagloo} generalized the result to all $\alpha\not\in\frac{1}{2}+\mathbb{Z}$ by using ideas of Nishioka \cite{Nishioka} in his proof of strict disintegratedness of $P_I$\footnote{{The definitions of geometric triviality and strict disintegratedness of a strongly minimal set $X$ over $\mathbb{C}(t)$ used in the rest of this paper are taken to be those given in the conclusions of Fact \ref{trivialityP2}, with $X_{II}(\alpha)$ replaced by $X$}.}. The central idea relies on showing that for $\alpha\not\in\frac{1}{2}+\mathbb{Z}$, the model-theoretic algebraic closure coincides with the definable closure, hence giving orthogonality to Manin Kernels:

\begin{fact}\label{Interdefinability}
Let $\alpha,\beta\in\mathbb{C}$ be such that $\alpha\not\in \frac{1}{2}+\mathbb{Z}$. Let $K$ be any differential field extension of $\mathbb{C}(t)$ and let $y\in X_{II}(\alpha)\setminus K^{alg}$, $z\in X_{II}(\beta)\setminus K^{alg}$. If $z\in K\gen{y}^{alg}$, then $z\in \mathbb{C}(t)^{alg}\gen{y}$. In particular, if $\beta\not\in \frac{1}{2}+\mathbb{Z}$, then $y$ and $z$ interalgebraic over $K$ implies that they are interdefinable over $\mathbb{C}(t)^{alg}$.
\end{fact}
We have stated here a more general version of what is proved in \cite[Proposition 3.3]{Nagloo}. However, it is not hard to see that $\alpha=\beta$ is never used in the proof. On the other hand we also have to look at the half-integers.
\begin{remark}\label{mistake}
We point out a mistake in the proof of geometric triviality for $X_{II}(\alpha)$ when $\alpha\in\frac12+\mathbb Z$ in \cite{FreitagNagloo}. In that proof, strong minimality of $X_{II}(\alpha)$ is used to conclude that if  $y\in X_{II}(\alpha)$, then a nonconstant polynomial $p\in K[y,y']$ cannot divide its derivative $p'$. However, as we have just seen, this is false. Thus, removing the exceptional subvariety gives strong minimality, but does not eliminate the corresponding invariant polynomial. Consequently, the argument in \cite{FreitagNagloo} does not establish geometric triviality at the half-integer parameters. We will give an alternate proof in Section 3. On the other hand, that the invariant polynomial exists will play a crucial role in the proof of Theorem \ref{halfIntegersTheorem}.
\end{remark}

Finally we state the result of Freitag and Nagloo already mentioned in the introduction. For this we will need the notion of nonorthogonality.

\begin{defn}
    Let $X$ and $Y$ be two strongly minimal sets. We say that $X$ and $Y$ are \emph{nonorthogonal} if there is a definable relation $R\subset X\times Y$ such that the projection of $R$ to $X$ is cofinite, the projection of $R$ to $Y$ is cofinite, and both of these projection maps are finite-to-finite.  
\end{defn}

We will refer to such an $R$ as a definable finite-to-finite relation. If $X$ and $Y$ are both defined over $K$ and there is no $K$-definable finite-to-finite relation $R\subset X\times Y$, we say that $X$ and $Y$ are \emph{weakly orthogonal}. In the case of geometrically trivial strongly minimal types, nonorthogonality and weakly nonorthogonality coincide (see for example \cite[Corollary 5.5,
chapter 2, section 5]{Pillay-book}).

\begin{fact}
    Let $X$ and $Y$ be two strongly minimal, geometrically trivial types defined over $K$. Then $X$ and $Y$ are nonorthogonal if and only if they are weakly nonorthogonal. 
\end{fact}

By taking $X=Y$ in the above we see that indeed in the conclusions of Fact \ref{trivialityP2}, $K$ can be replaced by $\mathbb{C}(t)$. We also note that if $X$ and $Y$ are two nonorthogonal strongly minimal sets, then $X$ is geometrically trivial if and only if $Y$ is geometrically trivial. As noted in the introduction, we have from \cite{FreitagNagloo}:

\begin{fact}
Suppose that $\alpha, \beta \in \mathbb{C}$ with $\alpha,\beta$ transcendental. Then $P_{II} ( \alpha )$ is nonorthogonal to $P_{II} (\beta)$ if and only if $\beta - \alpha \in \mathbb{Z}$ or $\beta + \alpha \in \mathbb{Z}$. In the case that $\beta - \alpha \in \mathbb{Z}$ or $\beta + \alpha \in \mathbb{Z}$, the only algebraic relations between solutions of $P_{II} ( \alpha)$ and $P_{II} (\beta)$ are induced by the classical B\"acklund transformations. 
\end{fact}
Finally, our proof will use a strong form of the Painlev\'e property, the proof of which can be found in \cite[Theorem 2.1, Theorem 18.1]{GLS}
\begin{fact}\label{strongPainleveProperty}Let $\alpha\in\mathbb{C}$.
\begin{enumerate}
\item All local solutions of $P_{II}(\alpha)$ can be analytically continued to single-valued meromorphic solutions in the complex plane.
\item Every pole of a solution of $P_{II}(\alpha)$ is simple and has residue $\pm 1$.
\end{enumerate}
\end{fact}

\section{The half integers and the Gambier construction.}
We now aim to show how geometric triviality of $X_{II}(\alpha)$ at $\alpha\in\frac{1}{2}+\mathbb{Z}$ follows from that at $\alpha\in\mathbb{Z}$. The transformation and construction below are already well-known (cf. \cite{ClarksonJoshiPickering1999}). Using the B\"acklund transformation we only need to focus on $X_{II}(\frac{1}{2})$. Recall that we have the B\"acklund transformation $S:X_{II}(0)\rightarrow X_{II}(0)$ with $S(y):=-y$. We write down the differential equation expressing the quotient $X_{II}(0)/\gen{S}$ by setting

\[
u=-\frac{y'}y,\qquad w=y^2.
\]

Using $y''=2y^3+ty$ we get
\[
\begin{aligned}
u'
&=-\frac{y''}{y}+\left(\frac{y'}y\right)^2\\
&=u^2-t-2w,
\end{aligned}
\]
as well as\[
w'=2yy'=-2uw.
\]

Thus the differential system for $X_{II}(0)/\gen{S}$ is given as
\[
\begin{cases}
u'=u^2-t-2w\\
w'=-2uw\\
w\neq0
\end{cases}
\tag{$G$}
\] 
Let $Y$ be the set defined in $\mathcal{U}$ by (G). Then by construction $Y$ is strongly minimal and geometrically trivial, since $X_{II}(0)$ is.

On the other hand, we can also consider the system
\[
\begin{cases}
u'=u^2-t-2w\\
w'=-2uw\\
\end{cases}
\tag{$G^*$}
\]
on its own right. It has an invariant curve $w=0$ corresponding to the Riccati fiber $u'=u^2-t$. It is also in $K$-definable bijection with the differential equation
\begin{equation}\label{GambierODE}
u''=2u^3-2tu-1.
\end{equation}
as can easily be seen after eliminating $w$. By construction, the Riccati equation $u'=u^2-t$ is an order $1$ subvariety of equation (\ref{GambierODE}). In particular, the set $Y^{*}$ defined by 
\[u''=2u^3-2tu-1\qquad\text{ and }\qquad u'\neq u^2-t\] is in $K$-definable bijection with $Y$ and hence is strongly minimal and geometrically trivial. 

It turns out that a change of derivation and variables allows us to rewrite equation (\ref{GambierODE}) as $P_{II}(\frac{1}{2})$.
\begin{prop}
$X_{II}(\frac{1}{2})$, and consequently $X_{II}(\alpha)$ for all $\alpha\in\frac{1}{2}+\mathbb{Z}$, is geometrically trivial.
\end{prop}
\begin{proof}
Let $u$ with $u''=2u^3-2tu-1$ be as above. Choose $c\in\mathbb C$ such that $c^3=-\frac{1}{2}$, and change the derivation on $\mathcal{U}$ by setting $\partial=cD$. Let
\[
s=\frac{t}{c}\qquad \text{and}\qquad z=cu.
\]
Then $\partial s=1$, and
\[
\begin{aligned}
\partial^2z
&=c^3D^2u\\
&=c^3(2u^3-2tu-1)\\
&=2z^3+sz+\frac12.
\end{aligned}
\]

Furthermore, it is not hard to see that under this change of derivation and variables, the Riccati equation $u'=u^2-t$ identifies with the Riccati equation $\partial z=z^2+\frac{s}{2}$. Since $Y^{*}$ is strongly minimal and geometrically trivial and we have only taken a constant change of derivation and variables we obtain the desired result.
\end{proof}

\section{Puiseux series solutions at $t=\infty$}

We begin by isolating all the Puiseux series solutions of $P_{II}(\alpha)$ at $t=\infty$. Similar formal expansions can be found in \cite{JoshiMazzocco2003} and are related to the classical work of Boutroux \cite{Boutroux1913}. For now, we will not restrict to $P^{*}_{II}$, that is, we allow the solutions to possibly be non generic over $K$. Let us fix some notation for fields of formal power series.

\begin{notn} In this paper, we will use the following (mostly) standard notation:
\begin{enumerate}
\item We write $\mathbb{C}[[\tau]]$ for the ring of formal power series in variable $\tau$, and $\mathbb{C}((\tau))$ for the field of formal Laurent series.
\item We use the big-$O$ notation for formal Laurent series. More precisely, the notation $f(\tau)=O(\tau^r)$ means that every exponent occurring in $f$ is at least $r$.
\item We write $\mathcal{P}_{\infty}=\bigcup_{n\geq1}\mathbb C((t^{-1/n}))$ for the field of formal Puiseux series at $t=\infty$. By (2), if $f\in \mathcal{P}_{\infty}$, then $f(t)=O(t^r)$ means that every exponent occurring in $f$ is at most $r$.
\end{enumerate}
\end{notn}

The first important observation:
\begin{lemma}
\label{Puiseuxt=inftyorder}
    Every $y\in \mathcal{P}_{\infty}\setminus\{0\}$ which is a solution of $P_{II}(\alpha)$ has order $\frac{1}{2}$ or $-1$.
\end{lemma}
\begin{proof}
    Let $y$ be a solution of $P_{II}(\alpha)$ and assume $y\in \mathcal{P}_{\infty}\setminus\{0\}$. We write $y=\sum\limits_{i=0}^{\infty}a_it^{e-\frac{i}{n}}$ where $a_{0}\neq 0$, $n\geq 1$ and $e\in\mathbb{Q}$ is the order, i.e., $e:=\text{order}_{\infty}(y)$. Since $y$ is a solution of $y''=2y^3+ty+\alpha$, we have
    \begin{equation*}
        \sum\limits_{i=0}^\infty \text{\tiny$(e-\frac{i}{n})\left((e-1)-\frac{i}{n}\right)$}a_it^{(e-2)-\frac{i}{n}}=\sum\limits_{i,j,k}{\tiny\text{$2a_ia_ja_k$}}t^{3e-\frac{(i+j+k)}{n}}+\sum\limits_{i=0}^\infty a_it^{(e+1)-\frac{i}{n}}+\alpha.
    \end{equation*}
  \noindent On the right-hand side, the orders of the first two summands are $3e$ and $e+1$, respectively, while the third has order $0$ if $\alpha\neq 0$. 
  
  If $e=0$, then the order of the left-hand side, provided it is nonzero, is less than $-2$. This is a contradiction, since the highest power on the right-hand side is $1$. 
  
  If $e<-1$, then $e+1<0$ and $3e<0$. If $\alpha\neq 0$, the highest power on the right-hand side is $0$, while the order of the left-hand side is at most $e-2<0$, a contradiction. If $\alpha=0$, the highest power on the right-hand side is $e+1$, and comparing orders gives $e-2=e+1$, again a contradiction.

If $-1<e<0$, then $e+1>0$ and $3e<0$. Hence the highest power on the right-hand side is $e+1$, while the order of the left-hand side is at most $e-2$, a contradiction. If $e=-1$, then the order is as required.

If $e>0$, the highest power on the right-hand side is either $3e$ or $e+1$. If $3e\neq e+1$, then this highest power is positive, while the order of the left-hand side is at most $e-2$, a contradiction. Thus $3e=e+1$, giving $e=\frac{1}{2}$.

Hence, the order of of $y$ as Puiseux series at $t=\infty$ is $\frac{1}{2}$ or $-1$. 
\end{proof}

{We now aim to find all the Puiseux series solutions of $P_{II}(\alpha)$ at $t=\infty$. It will be more convenient for a moment to perform a change of variable $x$, and work with ``solutions at $x=0$''. For any given Puiseux series $y(t)$ of $P_{II}(\alpha)$ at $t=\infty$, we can let $x=t^{-\frac{3}{2}}$, and define $u(x)=x^{1/3}y(x^{-2/3})$. Then $u(x)$ is a formal Puiseux series at $x=0$. Note that since $y$ has order $\frac{1}{2}$ or $-1$, we have that $u$ has $0$ or $1$ as its least power of $x$. In particular, for any Puiseux series solution of $P_{II}(\alpha)$ at $t=\infty$, the corresponding Puiseux series at $x=0$ (as defined above) is in the ring of Puiseux series $\bigcup_{m\geq 1}\mathbb{C}[[x^{\frac{1}{m}}]]$.}

It is not hard to see that $u(x)$ satisfy the following differential equation
\begin{equation}
\label{eqn:Puiseuxatx=0}
    2u^3+u+\alpha x+\frac{1}{4}x^2u-\frac{9}{4}x^3u_x-\frac{9}{4}x^4u_{xx}=0.
\end{equation}

\noindent Indeed differentiating $y=t^{1/2}u$ and using $x'=-\frac{3}{2}t^{-\frac{5}{2}}$ gives
\begin{equation*}
  y'=  \frac{1}{2}t^{-\frac{1}{2}}u-\frac{3}{2}t^{-2}u_x
\end{equation*} and 
\begin{equation*}
    y''=-\frac{1}{4}t^{-\frac{3}{2}}u+\frac{9}{4}t^{-3}u_x+\frac{9}{4}t^{-\frac{9}{2}}u_{xx}.
\end{equation*}

\noindent As $y$ is a solution to $P_{II}(\alpha)$, we have that \[-\frac{1}{4}t^{-\frac{3}{2}}u+\frac{9}{4}t^{-3}u_x+\frac{9}{4}t^{-\frac{9}{2}}u_{xx}=2t^{\frac{3}{2}}u^3+t^{3/2}u+\alpha.\]

\noindent Multiplying the expression by $x=t^{-\frac{3}{2}}$ gives the desired equation (\ref{eqn:Puiseuxatx=0}). When $x=0$, the constant term of $u(x)$, namely c, must satisfy $2c^3+c=0$. Hence $c=0$ or $c^2=-\frac{1}{2}$. We call these 3 solutions the \emph{potential leading coefficients} of Puiseux series solutions at $t=\infty$, and denote them by $c_0:=0$, $c_{+}=\frac{i}{\sqrt{2}}$ and $c_{-}=-\frac{i}{\sqrt{2}}$. 

Recall that using Lemma \ref{Puiseuxt=inftyorder}, we can restrict to the ring of Puiseux series $\bigcup_{m\geq 1}\mathbb{C}[[x^{\frac{1}{m}}]]$.
\begin{lemma}\label{existencePuiseux}
 For every $c$, a potential leading coefficient of Puiseux series solutions at $t=\infty$, there is a unique  Puiseux series $u(x)\in\bigcup_{m\geq 1}\mathbb{C}[[x^{\frac{1}{m}}]]$ satisfying (\ref{eqn:Puiseuxatx=0}) with $c$ as its constant term. Moreover, $u(x)\in\mathbb{C}[[x]]$.
\end{lemma}
\begin{proof}
Let $c$ be a potential leading coefficient and let $A(U)=2U^3+U.$ It is not hard to see that at the three potential leading coefficients satisfy, $A'(c)\neq0$. Indeed we have
\[
A'(c_0)=1,
\qquad
A'(c_\pm)=6c_\pm^2+1=-2.
\]
This will play an important role in our proof. Assume $u=c+\sum\limits_{n=1}^\infty a_nx^{\frac{n}{m}}$ for $a_i\in\mathbb{C}$ and some denominator $m\geq1$. Let us write $w:=\sum\limits_{n=1}^\infty a_nx^{\frac{n}{m}}$, so that $u=c+w$, and since $A(c)=0$ we get that
\[A(u)=6c^2w+6cw^2+2w^3+w.\]
So if $c=c_0$, then $A(u)=w+2w^3$. While if $c=c_{\pm}$, then $A(u)=-2w+6cw^2+2w^3$. In any case we get
\begin{equation}\label{eqn:A(U)}A(u)=A'(c)w+6cw^2+2w^3.\end{equation}
Let us write $\mathcal{L}$ for the differential operator
\[\mathcal{L}:=\frac{1}{4}x^2-\frac{9}{4}x^3\frac{d}{dx}-\frac{9}{4}x^4\frac{d^2}{dx^2},\]
so that $\mathcal{L}(u)+A(u)+\alpha x=0$ using equation (\ref{eqn:Puiseuxatx=0}). For any $r\in\mathbb{Q}$, with $r>0$ we get
\begin{equation}\label{operatorindex}\mathcal{L}(x^r)=\left(\frac{1}{4}-\frac{9}{4}r-\frac{9}{4}r(r-1)\right)x^{r+2}=\frac{1-9r^2}{4}x^{r+2}.\end{equation}
In other words, the differential operator $\mathcal{L}$ raises every exponent by 2. Also notice that $\mathcal{L}(c)=\frac{c}{4}x^2$.  

We now study the coefficients of $x^{\frac{n}{m}}$ in each term of $\mathcal{L}(u)+A(u)+\alpha x$. Using equation (\ref{eqn:A(U)}), the contribution from $A(u)$ is
\[A'(c)a_n+6c\sum_{\substack{i+j=n \\ i,j\geq 1}}a_ia_j+2\sum_{\substack{i+j+k=n \\ i,j,k\geq 1}}a_ia_ja_k.\]
On the other hand, using equation (\ref{operatorindex}) (with $\frac{n}{m}=r+2$), the contribution from $\mathcal{L}(u)$ is 
\[\frac{1}{4}\left(1-9\left(\frac{n-2m}{m}\right )^2\right)a_{n-2m},\]
where here we set $a_k=0$ if $k<0$. Notice also that if $n=2m$ (that is $r=0$), then the coefficient of $x^2$ in $\mathcal{L}(u)$ is $c/4$ as mentioned above. Finally the contribution of $\alpha x$ to the coefficients of $x^{\frac{n}{m}}$ is $\alpha\delta_{m,n}$, where $\delta_{m,n}$ is the Kronecker delta symbol.

Since $\mathcal{L}(u)+A(u)+\alpha x=0$, we have that the coefficient of $x^{\frac{n}{m}}$ satisfies
\begin{equation}\label{eq:an}
\begin{aligned}
a_n=-\frac{1}{A'(c)}\Bigg(
&6c\sum_{\substack{i+j=n \\ i,j\geq 1}}a_ia_j
+2\sum_{\substack{i+j+k=n \\ i,j,k\geq 1}}a_ia_ja_k
+\alpha\delta_{m,n}\\
&+\frac{1}{4}\left(1-9\left(\frac{n-2m}{m}\right)^2\right)a_{n-2m}+c\delta_{2m,n}\Bigg).
\end{aligned}
\end{equation}

So for fixed $m\geq 1$, we see that $a_n$ is determined by the value of $a_1,\ldots,a_{n-1}$. Since $a_1=-\frac{\alpha\delta_{m,1}}{A'(c)}$ is uniquely determined, we have that each coefficient $a_n$ is uniquely determined, proving uniqueness of such Puiseux series for any fixed denominator $m$. 

Taking $m=1$ gives the existence of a solution $u_1\in\mathbb{C}[[x]]$ of (\ref{eqn:Puiseuxatx=0}) with constant term $c$. Since $u_1\in\mathbb{C}[[x^{1/m}]]$, for every $m\geq1$, uniqueness for each fixed $m$ shows that $u_1$ is the unique formal Puiseux solution with constant term $c$.



\end{proof}

In the above lemma, equation (\ref{eq:an}) gives an expression for $a_n$, the coefficient of $x^{n}$ in $u$. Recall that $y(x)=x^{-\frac{1}{3}}u(x)$, $x=t^{-\frac{3}{2}}$ and $u(x)=c+\sum\limits_{n=1}^\infty a_nx^n\in \mathbb{C}[[x]]$. Hence $$y(t)=ct^\frac{1}{2}+\sum\limits_{n=1}^\infty a_nt^{\frac{1-3n}{2}}$$ for $c=0$ or $c=c_{\pm}$. 

As in the proof of Lemma \ref{existencePuiseux}, let $A(U)=2U^3+U$ so that $A(c)=0$ while $A'(c)\neq 0$. When $c=0$, we have $A'(0)=1$ and the order of the Puiseux series solution at $t=\infty $ is $-1$. Computing the coefficients from equation (\ref{eq:an}), we obtain that $a_1=-\alpha$ and $a_2=0$. Using these coefficients in $a_n$, we obtain that $a_3=2\alpha(\alpha^2-1)$, $a_4=0$, and $a_5=-4\alpha(\alpha^2-1)(3\alpha^2-10)$. Hence, the unique Puiseux series solution at $t=\infty$ of $P_{II}(\alpha)$ with order $-1$ is 
\begin{equation}\label{eqn:uniquey0}
y_{\alpha,0}=-\alpha t^{-1}+2\alpha(\alpha^2-1)t^{-4}+O(t^{-7}).
\end{equation}

When $c=c_{\pm}$, we have $A'(c_\pm)=-2$, and the order of the Puisseux series solution at $t=\infty $ is $\frac{1}{2}$. Again, using equation (\ref{eq:an}) to compute coefficients, we have that $a_1=\frac{\alpha}{2}$ and $a_2=\frac{c_{\pm}(6\alpha^2+1)}{8}$. Hence, the two Puiseux series solutions at $t=\infty$ of order $\frac{1}{2}$ are 
\begin{equation}\label{eqn:uniqueycneq0}
    y_{\alpha,\pm}=c_{\pm}t^{\frac{1}{2}}+\frac{\alpha}{2}t^{-1}+\frac{c_{\pm}(6\alpha^2+1)}{8}t^{-\frac{5}{2}}+O(t^{-4}).
\end{equation}

In particular, there is no term of order strictly between $t^{\frac{1}{2}}$ and $t^{-1}$ for any Puiseux series solution at $t=\infty$ of $P_{II}(\alpha)$.

We call $y_{\alpha,+}$ and $y_{\alpha,-}$, the solutions corresponding to the potential leading coefficients $c_{+}$ and $c_-$ the large branches, and $y_{\alpha,0}$, the solution corresponding to $c=0$, the small branch.

\begin{rem}\label{remarkBranches} The following observation holds for the large and short branch solutions:
\begin{enumerate}
\item If $\alpha\not\in\frac{1}{2}\mathbb{Z}$, all solutions (so including $y_{\alpha,0}$ and $y_{\alpha,\pm}$) are generic over $K$.
\item If $\alpha\in\mathbb{Z}$, then the short branch solution is the unique rational solution $y_\alpha$ of $P_{II}(\alpha)$. Indeed, it is not hard to see that at $\alpha=0$, we have $y_{0,0}=0$ (by Lemma \ref{Puiseuxt=inftyorder} the leading order of a nonzero solution has to be $\frac{1}{2}$ or $-1$). The classical B\"acklund transformation $T_+$ and $T_-$ preserve the rational solutions. On the other hand, the large branch solutions are generic over $K$ (since all other solutions are).
\item If $\alpha\in\frac{1}{2}+\mathbb{Z}$, then the short branch solution $y_{\alpha,0}$ does not satisfy the Riccati locus $R_\alpha$ and hence $y_{\alpha,0}$ is generic over $K$. Indeed, at $\alpha=\frac{1}{2}$ the Riccati equation defining $R_{\frac{1}{2}}$ is $y'=y^2+\frac{t}{2}$. Substituting $y_{\alpha,0}$ in that equation gives a contradiction by comparing the leading powers of $t$ ($t^{-2}$ on the LHS and $t^1$ on the RHS). The classical B\"acklund transformation $T_+$ and $T_-$ preserve the short branch solutions. The long branch solutions, on the other hand, belong to $R_\alpha$ (this is a classical fact that the Riccati equation has such solutions cf. \cite[Theorem 3.3]{Clarkson2016}).
\item The change of derivation and variables used in the Gambier construction in Section 3 allows us to determine the branch solutions of the Riccati equation $u'=u^2-t$. Indeed, starting with the two long branches $z_{\frac{1}{2},\pm}(s)$ of $\partial z=z^2+\frac{s}{2}$, we get that  $u'=u^2-t$ has only two large branch solutions $u_\pm$ satisfying \[
u_\pm=\pm t^{1/2}+O(t^{-1}).
\]
This follows since $u_\pm=\frac{z_{1/2,\pm}}{c}
=\frac{c_\pm}{c^{3/2}}t^{1/2}+O(t^{-1})$ and $
\left(\frac{c_\pm}{c^{3/2}}\right)^2
=1.
$
\end{enumerate}    
\end{rem}
\begin{prop}
Let $\alpha\in\mathbb{C}$. Then $X_{II}(\alpha)$ is $\omega$-categorical: if $K$ is any differential field extending $\mathbb{C}(t)$ and $y_1,...,y_n\in X_{II}(\alpha)\setminus K^{alg}$ are distinct, then $K(y_1,y_1',....,y_n,y'_n)^{alg}$ contains only finitely many elements of $X_{II}(\alpha)$.  
\end{prop}
\begin{proof}
If not since $X_{II}(\alpha)$ is geometrically trivial, we get that for any solution $y\in X_{II}(\alpha)$, the field $\mathbb{C}(t)(y,y')^{alg}$ contains infinitely many distinct elements of $X_{II}(\alpha)$. So in particular that ``there exists three more distinct solutions of $X_{II}(\alpha)$ in $\mathbb{C}(t)(y,y')^{alg}$'' is expressible by a first order formula $\sigma(u)$ over $\mathbb{C}(t)$ and $\sigma(u)$ is in the type of $y$ over $\mathbb{C}(t)$. In particular, we have that if $\sigma(z)$ is true in $\mathcal{U}$, then $z$ is a $\mathbb{C}(t)$-generic solution of $X_{II}(\alpha)$ and there exists three more distinct solutions of $X_{II}(\alpha)$ in $\mathbb{C}(t)(z,z')^{alg}$.

Since $X_{II}(\alpha)$ is strongly minimal and all solutions are generic, it follows that $\sigma(\hat y)$ is true in $\mathcal{U}$, where $\hat y$ is a generic branch solution in $X_{II}(\alpha)$. So there are three more distinct solutions of $X_{II}(\alpha)$ in $\mathbb{C}(t)(\hat y,\hat y')^{alg}\subseteq\mathcal{P}_{\infty}$ (using that the latter is algebraically closed). This contradicts that there are only 3 Puiseux series solutions of $P_{II}(\alpha)$ at $t=\infty$.
\end{proof}
\section{Algebraic relations within an equation}
Next we show that the only self-maps of $X_{II}(\alpha)$ over $K=\mathbb{C}(t)^{alg}$ occurs at $\alpha\in\mathbb{Z}$ and are given by the classical B\"acklund transformations\footnote{Note that $T_+=T_-^{-1}$.}
\[
S_{\alpha}:=T_+^{\alpha}\circ S\circ T_-^{\alpha}.
\] We begin by showing that for $\alpha,\beta\not\in\frac{1}{2}\mathbb{Z}$, the $K$-definable bijections have a special form.
\subsection{The Nishioka normal form of definable permutations of $P_{II}$}\label{section:NishiokaonPII}
Following the ideas of Nishioka in \cite{Nishioka}, we show that definable bijections in the $P_{II}$ family have a very special form, which we call the {\em Nishioka normal form}. Throughout we keep the notation $K=\mathbb{C}(t)^{alg}$ and let $L$ be any differential field extension of $\mathbb{C}(t)$ with $L\subseteq K$. 

Let us quickly recall two essential tools introduced by Nishioka to approach this problem. We assume throughout that $\alpha\not\in \frac{1}{2}+\mathbb{Z}$, since Condition (J) (Fact \ref{SM-Invariant}) is needed as part of our assumption. We direct the reader to Nishioka's original work in \cite{Nishioka} for more details.

For a given $y\in X_{II}(\alpha)$, we define a weight function and valuation as follows. The weight function $w$ on $L[y,y']$ is given by the following conditions and properties:

\begin{enumerate}
    \item $w(y)=1$, $w(y')=2$, $w(u)=0$ for all $u\in L\setminus\{0\}$,
    \item $w(P+Q)\leq\max\{w(P),w(Q)\}$ where the equality holds if $w(P)\neq w(Q)$.
    \item $w(PQ)=w(P)+w(Q).$
    \item $w(P')\leq w(P)+1$.
\end{enumerate}

Given $p\in L[y,y']$ irreducible, we will also consider the discrete valuation $\upsilon_p:L(y,y')\to\mathbb{Z}\cup\{\infty\}$ defined by $\upsilon_p(f)=m$ where $f=p^m\frac{h}{g}$ for some $h,g\in L[y,y']$ with $p\nmid h$ and $p\nmid g$ (i.e. $\upsilon_p(h)=\upsilon_p(g)=0$). If $f=0$, we set $\upsilon_p(0)=\infty$.  As the valuation is discrete, we have the properties: 
\begin{enumerate}
    \item $\upsilon_p(gh)=\upsilon_p(g)+\upsilon_p(h)$, and
    \item $\upsilon_p(g+h)\geq \min\{\upsilon_p(g),\upsilon_p(h)\}.$
\end{enumerate}

\begin{thm}\label{NishiokaNormalForm}
Let $\alpha,\beta\not\in \frac{1}{2}+\mathbb{Z}$ and let $f:X_{II}(\alpha)\rightarrow X_{II}(\beta)$ be a $L$-definable  bijection. Then $f$ has the form
\begin{equation}\label{eq:NishiokaNormalForm}
f(Y,Y')=\varepsilon Y+b+\frac{q(Y,Y')'}{q(Y,Y')},
\end{equation}
where $\varepsilon\in\{\pm1\}$, $b\in L$ and $q\in L(y,y')$.
\end{thm}

We will sometimes write $\text{dlog}\;q:=\frac{q'}{q}$ for the logarithmic derivative of $q$. We first need several lemmas.

\begin{lemma}\label{NoPolynomial}
Let $\alpha,\beta\not\in \frac{1}{2}+\mathbb{Z}$  and assume $f:X_{II}(\alpha)\rightarrow X_{II}(\beta)$ is a $L$-definable bijection. If $f\in L[Y,Y']$, then either $f=\text{id}$ or $f=S$, where $S(Y)=-Y$ is the classical B\"acklund transformation $X_{II}(\alpha)\rightarrow X_{II}(-\alpha)$.
\end{lemma}
\begin{proof}
Let $f\in L[Y,Y']$ be such that $f:X_{II}(\alpha)\rightarrow X_{II}(\beta)$. Let $y\in X_{II}(\alpha)$ and consider the weight function $w$ on $L[y,y']$ as defined above. Since $z:=f(y,y')\in L[y,y']$, we can consider the weights on both sides of $z''=2z^3+tz+\beta$. Doing so gives $w(z)\leq 1$ and hence $z=uy+v$ with $u,v\in L$. A quick computation, substituting $z=uy+v$ in  $z''=2z^3+tz+\beta$ and using the fact that $y\in X_{II}(\alpha)$, yields $u=\pm 1$ and $v=0$ from which the result follows.
\end{proof}

\begin{lemma}\label{AlmostNormalForm}
Let $\alpha,\beta\not\in \frac{1}{2}+\mathbb{Z}$  and assume $f:X_{II}(\alpha)\rightarrow X_{II}(\beta)$ be a $L$-definable bijection. Then $f$ has the form
\[
aY+b+\text{dlog}\;q,
\]
where $a,b\in L$ and $q\in L(Y,Y')$.
\end{lemma}
\begin{proof}
Let $f:X_{II}(\alpha)\rightarrow X_{II}(\beta)$ be a $L$-definable bijection and let $y\in X_{II}(\alpha)$. Using Lemma \ref{NoPolynomial}, we may assume that $f\not\in L[y,y']$. Let us write $z:=f(y,y')=\frac{h(y,y')}{g(y,y')}$, where $g,h\in L[Y,Y']$ and $\gcd(g,h)=1$. So for some $p\in L[y,y']$, we have that $\upsilon_p(z)<0$, where $\upsilon_p$ is the discrete valuation on $L(y,y')$ defined above. Since for any such $p$, we have that $\upsilon_p(z'')=\upsilon_p(z)-2$, it follows that $\upsilon_p(z)=-1$. In particular, the denominator $g$ is squarefree. 

Since $z''=2z^3+tz+\beta$ we get that \[
g^2h''-hgg''-2gg'h'+2h(g')^2
 =2h^3+thg^2+\beta g^3,
\]
that is \[
g(gh''-hg''-2g'h'-thg-\beta g^2)
 =2h(h^2-g'^2).
\]
Hence, $g|(h-g')(h+g')$. Since $g$ is squarefree, it follows that each irreducible factor of $g$ divides exactly one of the two factors $h-g'$ or $h+g'$. In other words $g=g_1g_2$ for some $g_1,g_2\in L[Y,Y']$ with $(g_1,g_2)=1$ and such that $g_1|(h-g')$ and $g_2|(h+g')$. 

Since $g'=g_1'g_2+g_2'g_1$ we have that $g'\equiv g_1'g_2 \Mod{g_1}$ and $g'\equiv g_2'g_1 \Mod{g_2}$. So together with $g_1|(h-g')$ and $g_2|(h+g')$, it follows that $h\equiv g_1'g_2 \Mod{g_1}$ while  $h\equiv -g_2'g_1 \Mod{g_2}$. A simple computation gives that $g|h(h-(g_1'g_2-g_2'g_1))$ and since $(g,h)=1$, by setting $q:=\frac{g_1}{g_2}$, we get that $z=l+\text{dlog}\;q$ where $l\in L[y,y']$. It remains to show that $w(l)\leq 1$.

Let $d:=w(z)$ and consider the weights on both sides of $z''=2z^3+tz+\beta$. It is not hard to see that $d\leq 1$. Indeed, if $d>1$, then $w(2z^3)=3d$ is greater than the weight of any other terms on either side of the equation (we have $w(z'')\leq d+2$, $w(tz)\leq d$ and $w(\beta)=0$) which is a contradiction. Finally, since $z=l+\text{dlog}\;q$, we get $w(l)\leq 1$ and we are done.
\end{proof}

We are now ready for the proof of Theorem \ref{NishiokaNormalForm}. All we have to show is that the element $a\in L$, given in Lemma \ref{AlmostNormalForm}, satisfies $a=\pm 1$.
\begin{proof}[Proof of Theorem \ref{NishiokaNormalForm}]
Let $f:X_{II}(\alpha)\rightarrow X_{II}(\beta)$ be given and assume $y\in X_{II}(\alpha)$ and $z\in X_{II}(\beta)$ are such that $z=f(y)$. Since $y$ and $z$ are interdefinable over $L$, we also have $y=F(z)$ for some $L$-definable bijection satisfying $F\circ f=\text{id}$, and using Lemma \ref{AlmostNormalForm} we can write
\[
z=ay+b+\text{dlog}\;q\:\text{ and }\:y=a_1z+b_1+\text{dlog}\;r
\]
for some $a,b,a_1,b_1\in L$ and $q,r\in L(Y,Y')$. Recall from the proof of Lemma \ref{AlmostNormalForm} that $q=\frac{g_1}{g_2}$ where $g_1,g_2\in L[y,y']$ are squarefree. The same is true of $r$. After substituting $z$ and $z'$ in $r$, we write the resulting rational function as $s\in L(y,y')$. After further substitution we get
\begin{equation}\label{eqnPM}
(1-aa_1)y=a_1b+b_1+a_1\text{dlog}\;q+\text{dlog}\;s.
\end{equation}
Notice that if $\text{dlog}\;q\in L$, then $z=ay+u$ with $u\in L$ and Lemma \ref{NoPolynomial} gives the result. So we can assume $\text{dlog}\;q\not\in L$ and hence in particular $q\not\in L$. Let $p\in L[y,y']\setminus L$ be irreducible such that $v_p(q)\neq 0$. Since $g_1$ and $g_2$ are relatively prime and squarefree, we get that $m:=\upsilon_p(q)=\pm 1$. Let $n:=\upsilon_p(s)$ and let us write 
\[q=p^mu\:\text{ and }s=p^nv\]
where $u,v\in L(y,y')$ are such that $p$ divides neither the numerator nor denominator of $u$ or $v$. Then
\[\frac{q'}q=m\frac{p'}p+\frac{u'}u,
\;\;\;
\frac{s'}s=n\frac{p'}p+\frac{v'}v
\]
and equation (\ref{eqnPM}) becomes
\[
(1-aa_1)y-a_1b-b_1
=
(a_1m+n)\frac{p'}p
+a_1\frac{u'}u+\frac{v'}v.
\]
Since $p$ does not divide the numerator nor the denominator of $u$ and $v$, we can clear the denominators by multiplying by $pT$ for some polynomial $T\in L[y,y']$ satisfying $p\nmid T$. But then reducing modulo $p$ gives $$(a_1m+n)p'T\equiv 0 \Mod{p}$$ and since we also have that $p\nmid p'$ (by Fact \ref{SM-Invariant}), we get that $a_1m+n=0$ (recall that $p\notin L$ and that $a_1m+n\in L$). Together with $m=\pm 1$, this also means that $a_1\in\mathbb{Z}$. The above argument also shows that for any irreducible polynomial $p\in L[y,y']\setminus L$, if the valuation $v_p(q)\neq 0$, then $a_1\upsilon_p(q)+\upsilon_p(s)=0$.

On the other end, let us assume that $p\in L[y,y']\setminus L$ is irreducible and such that $v_p(q)=0$. As before, we write $n:=\upsilon_p(s)$ and \[s=p^nv\]
where $v\in L(y,y')$ is such that $p$ does not divide its numerator and denominator (and same for $q$, since $v_p(q)=0$). This time equation (\ref{eqnPM}) becomes
\[
(1-aa_1)y-a_1b-b_1
=
n\frac{p'}{p}
+a_1\frac{q'}{q}+\frac{v'}v.
\]
As argued in the previous case, reducing modulo $p$ gives $$np'T\equiv 0 \Mod{p}$$ for some polynomial $T\in L[y,y']$ satisfying $p\nmid T$. Since $p\nmid p'$ (by Fact \ref{SM-Invariant}), we get $n=0$. So in all cases, we have established that if $p\in L[y,y']\setminus L$ is any irreducible polynomial, then $a_1\upsilon_p(q)+\upsilon_p(s)=0$ (and $a_1\in\mathbb{Z}$). In other words, $q^{a_1}s\in L$, so that $$a_1\text{dlog}\;q+\text{dlog}\;s=\frac{(q^{a_1}s)'}{(q^{a_1}s)}\in L.$$ Finally, equation (\ref{eqnPM}) gives $(1-aa_1)y\in L$ and hence $aa_1=1$, since $y\not\in L$. Since a symmetric argument also gives $a\in\mathbb{Z}$, we conclude that $a=\pm 1$ as required.
\end{proof}


\subsection{Definable permutations of $P_{II}$ on the branches} We are now ready to study self-maps. In Theorem \ref{NishiokaNormalForm}, we showed that every definable permutation is of Nishioka normal form. Using the branches in the previous section, we can show that the $b$ and $\text{dlog}(q)$ terms in these permutations are of order $-1$.  

\begin{lemma}\label{dlogisO-1}
Let $q\in K(Y,Y')$ with $q\not\equiv0$.
\begin{enumerate}
    \item If $\hat{y}\in\{y_{\alpha,0},y_{\alpha,+},y_{\alpha,-}\}$ is a branch solution with $q(\hat{y},\hat{y}')$ well defined and $q(\hat{y},\hat{y}')\neq 0$, then
    \[
    \text{dlog}(q(\hat{y},\hat{y}'))=O(t^{-1}).
    \]

    \item If $q\in\mathbb{C}(t)(Y,Y')$ and the short branch solution $y_{\alpha,0}$ is generic, then
    \[
    \text{dlog}(q(y_{\alpha,0},y_{\alpha,0}'))
    =\lambda t^{-1}+O(t^{-2})
    \]
    for some $\lambda\in\mathbb{Z}$.
\end{enumerate}
\end{lemma}

\begin{proof}
(1) Let $\hat{y}$ be one of the three branch solutions of $P_{II}(\alpha)$. By assumption we have that $q(\hat{y},\hat{y}')\neq 0$ and since $\mathcal{P}_{\infty}$ is a differential field extension of $K$, it follows that $Q(t)=q(\hat{y},\hat{y}')$ is a non-zero Puiseux series. Let $r\in\mathbb{Q}$ be the largest power of $t$ in $Q(t)$, so that $Q(t)=at^r(1+O(1))$. Differentiating gives, 
    \begin{align*}
        \text{dlog}(q)(\hat{y},\hat{y}')&=\frac{Q'(t)}{Q(t)}\\
        &= \frac{ar t^{r-1}(1+O(1))+at^r O(t^{-1})}{at^r(1+O(1))}\\
        &=\frac{r}{t}+O(t^{-1}).
    \end{align*}

\noindent (2) Let $\hat{y}=y_{\alpha,0}$ be a generic short branch solution. Notice by Lemma \ref{existencePuiseux} that $y_{\alpha,0}\in\mathbb{C}((t^{-1}))$\footnote{We use here that $u(x)\in\mathbb C[[x]]$, $x=t^{-3/2}$ and $y=t^{-1}u(x)$.} and that $\mathbb{C}((t^{-1}))$ is a differential field extension of $\mathbb{C}(t)$ (which can be easily checked). Hence, it follows that $Q(t):=q(\hat{y},\hat{y}')\in\mathbb{C}((t^{-1}))$ since $q$ is over $\mathbb{C}(t)$, and since $y_{\alpha,0}$ is generic $Q(t)\neq0$. In other words, 
$Q(t)=at^{\lambda}(1+O(t^{-1}))$ for $\lambda\in\mathbb{Z}$. Hence
\begin{align*}
   \text{dlog}(q)(\hat{y},\hat{y}')&= \frac{Q'(t)}{Q(t)}\\
    &=\frac{\lambda at^{\lambda-1}(1+O(t^{-1}))+at^{\lambda}O(t^{-2})}{at^{\lambda}(1+O(t^{-1}))}\\
    &=\lambda t^{-1}+O(t^{-2}),
\end{align*} and the result follows.
\end{proof}

\begin{lemma}\label{fixShort}
Suppose $\alpha,\beta\notin\frac{1}{2}\mathbb {Z}$, and $f:X_{II}(\alpha)\rightarrow X_{II}(\beta)$ is a $K$-definable bijection in Nishioka normal form $f(Y,Y')=\epsilon Y+b+\mathrm{dlog}(q)$. Then $$f(y_{\alpha,0})=y_{\beta,0},$$ where $y_{\alpha,0}$ and $y_{\beta,0}$ are the (generic) short branch solutions. Consequently, $b=O(t^{-1})$.
\end{lemma}
\begin{proof}
Since $f$ takes a Puiseux series solution to another, we have that $f$ gives a bijection between $\{y_{\alpha,0},y_{\alpha,+},y_{\alpha,-}\}$ and $\{y_{\beta,0},y_{\beta,+},y_{\beta,-}\}$. Using Lemma \ref{dlogisO-1}, we have that if $\hat y\in \{y_{\alpha,0},y_{\alpha,+},y_{\alpha,-}\}$, then $\text{dlog}(q(\hat{y},\hat{y}'))=O(t^{-1})$. So we can look at the difference 
\[
\begin{aligned}
f(y_{\alpha,+})-f(y_{\alpha,-})
&=\varepsilon(y_{\alpha,+}-y_{\alpha,-})+O(t^{-1})\\
&=2\varepsilon c\,t^{1/2}+O(t^{-1}),
\end{aligned}
\] which has leading coefficient $\pm 2c$, where $c^2=-1/2$. On the other hand, it is not hard to see that out of the leading terms of the differences $y_{\beta,i}-y_{\beta,j}$, with $i,j\in\{0,+,-\}$, only the difference of the long branches has leading coefficients $\pm 2c$. Therefore $f$ must map the pair of long branches onto the pair of long branches and $f(y_{\alpha,0})=y_{\beta,0}$.

Finally, since $b=y_{\beta,0}-\varepsilon y_{\alpha,0}-\text{dlog}(q(y_{\alpha,0},y'_{\alpha,0})$ and each term on the RHS is $O(t^{-1})$, we get that $b=O(t^{-1})$.
\end{proof}

\begin{thm}\label{disintegrated}
 The following holds:    
\begin{enumerate}
\item If $\alpha\not\in \frac{1}{2}\mathbb{Z}$, then $X_{II}(\alpha)$ is strictly disintegrated, that is the identity map is the unique $K$-definable permutation of $X_{II}(\alpha)$.
\item If $\alpha\in\mathbb{Z}$, then the only $K$-definable permutations of $X_{II}(\alpha)$ are the identity and $S_{\alpha}$.
\end{enumerate}  
\end{thm}
\begin{proof}
(1) If not, then by Fact \ref{trivialityP2} and \ref{Interdefinability} there is $f:X_{II}(\alpha)\rightarrow X_{II}(\alpha)$ is a $K$-definable permutation with $f\neq id$. By Lemma \ref{fixShort}, we have that $f(y_{\alpha,0})=y_{\alpha,0}$ where $y_{\alpha,0}$ is the short branch solution. But $y_{\alpha,0}$ is a generic solution of $X_{II}(\alpha)$. Hence $f(Y)=Y$, a contradiction.

\noindent (2) Assume $f:X_{II}(\alpha)\rightarrow X_{II}(\alpha)$ is a $K$-definable permutation in Nishioka normal form with $f\neq id$. Since $f\neq id$ and preserves Puiseux series solutions, we have that $f(y_{\alpha,+})=y_{\alpha,-}$, where $y_{\alpha,\pm}$ are the long branch solutions. But since $y_{\alpha,+}$ is generic, this completely determine $f$. So there is only one $K$-definable permutation of $X_{II}(\alpha)$ with $f\neq id$. So $f=S_{\alpha}$.
\end{proof}

\begin{cor}\label{halfIntegerDisintegrated}
Let $\alpha\in\frac{1}{2}+\mathbb {Z}$, then $X_{II}(\alpha)$ is strictly disintegrated
\end{cor}
\begin{proof}
This follows from the discussion in Section 3. The definable set $X_{II}(\frac{1}{2})$ inherits the same structure as that of $X_{II}(0)/\gen{S}$. By Theorem \ref{disintegrated}(2), it follows that $X_{II}(0)/\gen{S}$ is strictly disintegrated and we are done after applying the classical B\"acklund transformations to $X_{II}(\frac{1}{2})$.
\end{proof}

\section{Transformations between different equations}
We are now ready to prove our main theorem which we now recall
\begin{thm}
Let $L$ be any differential field extension of $\mathbb{C}(t)$. Suppose $y_1,\ldots,y_n$ are $L$-generic solutions of $X_{II}(\alpha_1),\ldots,X_{II}(\alpha_n)$, respectively, and $tr.deg_{L}L(y_1,y_1',\ldots,y_n,y_n')<2n$, then for some $i<j$, we have that $\alpha_i\pm\alpha_j\in\mathbb{Z}$, and $y_i=T(y_j)$ for a classical B\"acklund transformation $T$.   
\end{thm}

Let $y_1,\ldots,y_n$ be $L$-generic solutions of $P_{II}(\alpha_1),\ldots,P_{II}(\alpha_n)$ and assume $tr.deg_{L}L(y_1,y_1',\ldots,y_n,y_n')<2n$. Since $X_{II}(\alpha_1),\ldots,X_{II}(\alpha_n)$ are strongly minimal and geometrically trivial for some $i<j$ we have that $tr.deg_{\mathbb{C}(t)}\mathbb{C}(t)(y_i,y_i',y_j,y_j')=2$. In particular $X_{II}(\alpha_i)$ is nonorthogonal to $X_{II}(\alpha_j)$. So our theorem follows as soon as we classify the possible nonorthogonalities between pairs $X_{II}(\alpha)$ and $X_{II}(\beta)$. The remainder of the paper is devoted to this. We will break this problem into several cases.

\subsection{The easy cases} We start with some easy cases that readily follow from the analysis in Section 5.

\begin{prop}
Assume $\alpha,\beta\in\mathbb{Z}$ and $f:X_{II}(\alpha)\rightarrow X_{II}(\beta)$ is a definable bijection over $K$, then $f$ is a classical B\"acklund transformation.
\end{prop}
\begin{proof}
Suppose $\alpha,\beta\in\mathbb{Z}$ and $f:X_{II}(\alpha)\rightarrow X_{II}(\beta)$ is a definable bijection over $K$. Using the classical B\"acklund transformations, we may assume that $\alpha=0$. Recall that we already have a classical B\"acklund transformation $T:X_{II}(0)\rightarrow X_{II}(\beta)$. If $f=T$ then we are done. On the other hand, if $f\neq T$, then we can construct a non-trivial definable permutation $T^{-1}\circ f:X_{II}(0)\rightarrow X_{II}(0)$. Using Theorem \ref{disintegrated}, we get that $T^{-1}\circ f=S:=-y$. In other words $f=T\circ S$ is a classical B\"acklund transformation.
\end{proof} By Fact \ref{Interdefinability}, it follows that if $X_{II}(\alpha)$ is non-orthogonal to $X_{II}(\beta)$, when $\alpha,\beta\in\mathbb{Z}$, then non-orthogonality is witnessed using a classical B\"acklund transformation. Similarly, using strict disintegratedness, for $\alpha,\beta\in\frac{1}{2}+\mathbb{Z}$, nonorthogonality of $X_{II}(\alpha)$ and $X_{II}(\beta)$ is witnessed by a $K$-definable bijection $f:X_{II}(\alpha)\rightarrow X_{II}(\beta)$. We have
\begin{prop}
Assume $\alpha,\beta\in\frac{1}{2}+\mathbb{Z}$ and $f:X_{II}(\alpha)\rightarrow X_{II}(\beta)$ is a definable bijection over $K$, then $f$ is a classical B\"acklund transformation.
\end{prop}
\begin{proof}
Suppose $\alpha,\beta\in\frac{1}{2}+\mathbb{Z}$ and $f:X_{II}(\alpha)\rightarrow X_{II}(\beta)$ is a definable bijection over $K$. As before using the classical B\"acklund transformations, we may assume that $\alpha=\frac{1}{2}$ and we have a classical B\"acklund transformation $T:X_{II}(\frac{1}{2})\rightarrow X_{II}(\beta)$. We claim $f=T$. If not we get a non-trivial definable permutation $T^{-1}\circ f:X_{II}(\frac{1}{2})\rightarrow X_{II}(\frac{1}{2})$. But this contradicts Corollary \ref{halfIntegerDisintegrated}, and gives that $f=T$ is a classical B\"acklund transformation.
\end{proof}
It turns out that Fact \ref{Interdefinability} also reduces nonorthogonality of $X_{II}(\alpha)$ and $X_{II}(\beta)$, when $\alpha\not\in\frac{1}{2}\mathbb{Z}$ and $\beta\in\mathbb{Z}$, to the problem of classifying $K$-definable bijections. This time we get

\begin{prop}\label{parameterPreservation}
There are no $K$-definable bijections $f:X_{II}(\alpha)\rightarrow X_{II}(\beta)$ with $\alpha\not\in\frac{1}{2}\mathbb{Z}$ and $\beta\in\mathbb{Z}$.
\end{prop}
\begin{proof}
In this case we simply have to observe that if $y_{\alpha,i}\in\mathcal{P}_{\infty}$, with $i\in\{0,+,-\}$, is a solution in $X_{II}(\alpha)$ then $f(y_{\alpha,i})$ is again an element of $\mathcal{P}_{\infty}$ and a solution in  $X_{II}(\beta)$. So $f(y_{\alpha,i})$ is one of the branches $y_{\beta,i}$ with $i\in\{+,-\}$. It is hence impossible for $f$ to be a bijection.
\end{proof}So we have that $X_{II}(\alpha)$ and $X_{II}(\beta)$ are orthogonal whenever $\alpha\not\in\frac{1}{2}\mathbb{Z}$ and $\beta\in\mathbb{Z}$. The remaining easy case does not use Fact \ref{Interdefinability} instead that both definable sets are strictly disintegrated.
\begin{prop}\label{parameterPreservation}
Let $\alpha\not\in\frac{1}{2}\mathbb{Z}$ and $\beta\in\frac{1}{2}+\mathbb{Z}$. Then $X_{II}(\alpha)$ and $X_{II}(\beta)$ are orthogonal.
\end{prop}
\begin{proof}
Suppose $X_{II}(\alpha)$ and $X_{II}(\beta)$ are nonorthogonal. Since both sets are geometrically trivial, they are weakly nonorthogonal and so for any $y\in X_{II}(\alpha)$ there is $z\in X_{II}(\beta)$ so that $K(y,y')^{alg}=K(z,z')^{alg}$. Moreover, since $X_{II}(\alpha)$ and $X_{II}(\beta)$ are strictly disintegrated, it follows that $K(y,y')=K(z,z')$. In other words, we have a $K$-definable bijection between $X_{II}(\alpha)$ and $X_{II}(\beta)$. But $X_{II}(\alpha)$ has three Puiseux branch solutions while $X_{II}(\beta)$ has only one. As argued in the previous proof, this is impossible.
\end{proof}

\subsection{The general case} We now come to the heart of the proof of this section. Here we will need the full force of the Nishioka normal form. First we get
\begin{lemma}\label{unique}
 For fixed $\alpha,\beta\not\in\frac{1}{2}\mathbb{Z}$, there is at most one $K$-definable bijection from $X_{II}(\alpha)$ to $X_{II}(\beta)$.
\end{lemma}
\begin{proof}
Let $f:X_{II}(\alpha)\rightarrow X_{II}(\beta)$ be a $K$-definable bijection. From Lemma \ref{fixShort} we get $f(y_{\alpha,0})=y_{\beta,0}$, where $y_{\alpha,0}$ and $y_{\beta,0}$ are the the short branch solutions. Since $y_{\alpha, 0}$ is a generic solution of $X_{II}(\alpha)$, uniqueness of $f$ follows.
\end{proof}

\begin{lemma}\label{overCT}
Assume $\alpha,\beta\not\in\frac{1}{2}\mathbb{Z}$. Any $K$-definable bijection $X_{II}(\alpha)\rightarrow X_{II}(\beta)$ is defined over $\mathbb{C}(t)$.
\end{lemma}
\begin{proof}
This follows since for every $\sigma\in Gal(K/\mathbb{C}(t))$, we obtain by conjugating the coefficients another $K$-definable bijection $f^{\sigma}:X_{II}(\alpha)\rightarrow X_{II}(\beta)$. By Lemma \ref{unique}, we get that $f=f^{\sigma}$ for every $\sigma\in Gal(K/\mathbb{C}(t))$ and hence $f$ is defined over $\mathbb{C}(t)$.
\end{proof}
Our goal is to prove
\begin{thm}\label{mainTheorem}
Assume $\alpha,\beta\not\in\frac{1}{2}\mathbb{Z}$. Let $f:X_{II}(\alpha)\rightarrow X_{II}(\beta)$ be a $K$-definable bijection. Then $\beta\pm\alpha\in\mathbb{Z}$ and $f$ is generated by one of the classical B\"acklund transformations.
\end{thm}

Our strategy for the proof is straight forward. Since the parameters $\alpha$ and $\beta$ appear in the coefficients of the branch solutions $y_{\alpha,c}$ and $y_{\beta,c}$ with $c\in\{0,c_{\pm}\}$ we would like to compare the coefficients of the series in the expression
$$y_{\beta,c}=\varepsilon y_{\alpha,c}+b+\text{dlog}\;q$$
and hopefully obtain that $\beta\pm\alpha\in\mathbb{Z}$. 

It is unsurprising that we will now use the {\em Painlev\'e property} as part of our argument to classify all the possible transformations.

\begin{prop}\label{coefficientB}
Assume $\alpha,\beta\not\in\frac{1}{2}\mathbb{Z}$. Let $f:X_{II}(\alpha)\rightarrow X_{II}(\beta)$ be a $\mathbb{C}(t)$-definable bijection in Nishioka normal form
\[
f(Y)=\varepsilon Y+b+\text{dlog}\;q,
\]
where $b\in \mathbb{C}(t)$ and $q\in \mathbb{C}(t)(Y,Y')$. Then, if the Puiseux series of $b$ at $t=\infty$ is of the form $b(t)=b_{-1}t^{-1}+O(t^{-2})$, then $b_{-1}\in\mathbb{Z}$.
\end{prop}
\begin{proof}
By Fact \ref{strongPainleveProperty}, any local solutions $y$ of $P_{II}(\alpha)$ and $z$ of $P_{II}(\beta)$ extend meromorphically on the whole of $\mathbb{C}$ and has only simple
finite poles and integral residues. So for any $a\in\mathbb{C}$, we can write $y(t)=a_y(t-a)^{-1}+O(1)$ for $a_y\in\mathbb{Z}$, and similarly $z(t)=a_z(t-a)^{-1}+O(1)$ for $a_z\in\mathbb{Z}$. Assume that we are given two generic solutions such with $z=f(y)$. By construction, $$\text{dlog}\;q(y(t),y(t)')=\frac{Q'(t)}{Q(t)}$$
for some non-zero meromorphic function $Q$ on $\mathbb{C}$. If we write $Q(t)=(t-a)^mh(t)$ with $h(a)\neq0$ and $m\in\mathbb{Z}$, then using $Q'(t)=m(t-a)^{m-1}h(t)+(t-a)^mh'(t)$ we get from $b=z-\varepsilon y-\frac{Q'}{Q}$ that 
$$b=\frac{a_z-\varepsilon a_y-m}{t-a}+O(1).$$ In other words, every finite pole of $b$ is simple and has integer residue. Using the residue theorem for $\mathbb{P}^1$ (see \cite[Theorem 10.21]{Forster}), we have that $$b_{-1}:=\text{res}_{\infty}(b)=-\sum_{a\in\mathbb{C}}\text{res}_{a}(b)\in{\Large \mathbb{Z}}$$
as required.
\end{proof}

\begin{proof}[Proof of Theorem \ref{mainTheorem}]
Assume $\alpha,\beta\not\in\frac{1}{2}\mathbb{Z}$. Let $f:X_{II}(\alpha)\rightarrow X_{II}(\beta)$ be a $K$-definable bijection. By Lemma \ref{overCT}, we may assume $f$ is over $\mathbb{C}(t)$ and by Theorem \ref{NishiokaNormalForm} it has a Nishioka normal form with coefficients in $\mathbb{C}(t)$. Let $\delta=\beta-\varepsilon\alpha$, so that if $\varepsilon=1$ then $\delta=\beta-\alpha$ and if $\varepsilon=-1$ then $\delta=\beta+\alpha$. Using the short branch solutions, we have $$y_{\beta,0}=\varepsilon y_{\alpha,0}+b+\text{dlog}\;q(y_{\alpha,0},y'_{\alpha,0}).$$ By Lemma \ref{dlogisO-1}, $\text{dlog}\; q(y_{\alpha,0},y'_{\alpha,0})=\lambda t^{-1}+O(t^{-2})$ for some $\lambda\in\mathbb{Z}$. Also, from Proposition \ref{coefficientB}, we have that the Puiseux series of $b$ at $t=\infty$ is of the form $b(t)=b_{-1}t^{-1}+O(t^{-2})$ with $b_{-1}\in\mathbb{Z}$. So we get that $$-\beta=-\varepsilon\alpha+b_{-1}+\lambda$$ and hence $\delta:=\beta-\varepsilon\alpha\in\mathbb{Z}$.

Since when $\beta\pm\alpha\in\mathbb{Z}$ the classical B\"acklund transformation already exists, Lemma \ref{unique} gives the full result.
\end{proof}

\subsection{The special case of $X_{II}(0)$ and $X_{II}(\frac{1}{2})$} Finally, we look at the question of nonorthogonality of $X_{II}(\alpha)$ and $X_{II}(\beta)$ when $\alpha\in\mathbb{Z}$ and $\beta\in\frac{1}{2}+\mathbb{Z}$. Using the classical B\"acklund transformations we are done if we tackle the case of $X_{II}(0)$ and $X_{II}(\frac{1}{2})$. 

The existence of the Gambier construction hints at why this case is special: a change of variables carries some trace of the relation between $X_{II}(0)$ and $X_{II}(\frac{1}{2})$. Furthermore, since we do not have the Nishioka normal form for $X_{II}(\frac{1}{2})$, some of the previous techniques do not apply. The strategy in this case will be to do the Puiseux analysis at the level of the Riccati equations and then show that nonorthogonality of $X_{II}(0)$ and $X_{II}(\frac{1}{2})$ can be reduced to relation between these equations.

\begin{lemma}\label{riccati}
There are no $K$-definable bijection between the set defined by $u'=u^2-t$ and $z'=z^2+\frac{t}{2}$.
\end{lemma}
\begin{proof}
Let $X$ be the set defined by $u'=u^2-t$ and $Y$ that defined by $z'=z^2+\frac{t}{2}$. For contradiction, suppose $f:X\rightarrow Y$ is a $K$-definable bijection. Say $f\in K(U)$, so that for any solution $u\in X$ we have $f(u)'=f(u)^2+\frac{t}{2}$. We want to apply $f$ to the long branch solutions $u_\pm=\pm t^{1/2}+O(t^{-1})\in X$. As with the previous proof strategies, we would like to use the fact that on branch solutions $\text{dlog}=O(t^{-1})$. 

For $u\in X$, using $f(u)'=f_t(u)+(u^2-t)f_U(u)$, we can differentiate $f(u)'=f(u)^2+\frac{t}{2}$ with respect to $U$ and obtain

\begin{eqnarray*}
2f(u)f_{U}(u)&=&f_{tU}(u)+2uf_U(u)+(u^2-t)f_{UU}(u)\\
&=& f_{Ut}(u)+(u^2-t)f_{UU}(u)+2uf_U(u)\\
&=& f_u(u)'+2uf_U(u).
\end{eqnarray*}

Hence
\[
f(u)=u+\frac{1}{2}\frac{f_U(u)'}{f_U(u)}=u+\frac{1}{2}\text{dlog}\;f_U(u).
\]

{By Lemma \ref{dlogisO-1} on the large branch solutions $u_\pm$, we have that $\text{dlog}\;f_U(u_\pm)=O(t^{-1})$.} So $f(u_\pm)$ is a branch solution of $y'=y^2+\frac{t}{2}$ and $$f(u_\pm)=u_\pm+O(t^{-1}).$$ In other words the leading term of the branch solution $f(u_+)$ of $y'=y^2+\frac{t}{2}$ has leading coefficient $\pm 1$ contradicting that the branches of $y'=y^2+\frac{t}{2}$ have leading coefficient $c$ with $c^2=-1/2$ (see Remark \ref{remarkBranches}). 
 
\end{proof}
We are now ready to prove the last case.
\begin{thm}\label{halfIntegersTheorem}
Let $\alpha\in\mathbb{Z}$ and $\beta\in\frac{1}{2}+\mathbb{Z}$. Then $X_{II}(\alpha)$ and $X_{II}(\beta)$ are orthogonal.    
\end{thm}
\begin{proof}
    
Suppose not. Using the classical B\"acklund transformations we get that $X_{II}(0)$ and $X_{II}(\frac{1}{2})$ are (weakly) nonorthogonal. Say $R\subseteq X_{II}(0)\times X_{II}(\frac{1}{2})$ is the infinite finite-to-finite $K$-definable relation witnessing nonorthogonality. Since $X_{II}(\frac{1}{2})$ is strictly disintegrated for any $y\in X_{II}(0)$, there is one and only one $z\in X_{II}(\frac{1}{2})$ such that $(y,z)\in R$. In particular $z\in K(y,y')$. We now exploit the Gambier construction.
For any solution $y\in X_{II}(0)$ set
\[
u=-\frac{y'}{y}
\]
and rewrite $P_{II}(0)$ as the differential system
\[
y'=-uy,\qquad u'=u^2-t-2y^2.
\]
We work in $K(u)(y)=K(y,y')$, and for $z\in K(y,y')$ with $z\in X_{II}(\frac{1}{2})$, write
\[
z=\frac{P}{Q},
\qquad
P,Q\in K(u)[y],\quad \gcd(P,Q)=1.
\]

Since $z$ satisfies $P_{II}(1/2)$, we have $\left(\frac PQ\right)''=2\left(\frac PQ\right)^3+t\frac PQ+\frac{1}{2}$ and clearing denominators gives
\begin{equation}\label{equationPQ}
Q^2P''-PQQ''
-2QP'Q'+2PQ'^2
=2P^3+tPQ^2+\frac12Q^3.
\end{equation}

Notice that the ideal $(y)$ of $K(u)[y]$ is a differential ideal, since $Dy=-uy$ and $Du=u^2-t-2y^2$. Furthermore, it is not hard to see that $y$ divides neither $P$ nor $Q$. Indeed, if $y\mid Q$, then $y\mid Q'$, and equation (\ref{equationPQ}) implies $y\mid P^3$, contradicting $\gcd(P,Q)=1$. Similarly, if $y\mid P$, then $y\mid P'$ and $y\mid P''$. Equation (\ref{equationPQ}) therefore implies $y\mid Q^3$, again a contradiction.

Consequently, we may write
\[
P=p+yP_1,\qquad Q=q+yQ_1,
\qquad p,q\in K(u)^\times.
\]

Since $(y)$ is a differential ideal and $u'=u^2-t-2y^2$, we get a derivation $\delta$ on $K(u)[y]/(y)=K(u)$ satisfying 
\[\delta u=u^2-t.\] 
Furthermore, since
\[
\begin{aligned}
P'&\equiv\delta p\pmod{y},&
P''&\equiv\delta^2p\pmod{y},\\
Q'&\equiv\delta q\pmod{y},&
Q''&\equiv\delta^2q\pmod{y}.
\end{aligned}
\]

we can reduce equation (\ref{equationPQ}) modulo $(y)$ and get
\[
q^2\delta^2p-pq\delta^2q
-2q(\delta p)(\delta q)+2p(\delta q)^2
=2p^3+tpq^2+\frac12q^3.
\]

Putting $g=\frac pq\in K(u)^\times$ we see that, after dividing by $q^3$, the left hand side is simply $\delta^2g$ and so we get $\delta^2g=2g^3+tg+\frac{1}{2}$. Hence $g$ is a solution of $P_{II}(1/2)$ and since there are no algebraic solutions in $K$, we have that $g\notin K$. Furthermore, since $g,\delta g\in K(u)$ we have that $g$ is a nonalgebraic solution of $P_{II}(1/2)$ of order one. Hence
\[
\delta g=g^2+\frac{t}{2}.
\]
It is well known that any three distinct solutions of $\delta u=u^2-t$ are algebraically independent over $K$ (cf. \cite[Proposition 2.2]{Nag2}) and so $u\in K(g)$ (i.e., not just in $K(g)^{alg}$). It follows that $g$ is a $K$-definable bijection between the sets defined by $\delta u=u^2-t$ and $\delta g=g^2+\frac{t}{2}$. This contradicts Lemma \ref{riccati}. Hence $X_{II}(0)$ and $X_{II}(1/2)$ are orthogonal and the result follows.
\end{proof}

\end{document}